%% file: paper.tex
\documentclass[10pt,a4paper,twoside]{article}
\usepackage{comcon}
\usepackage{graphicx}

	\def\fpicps #1#2#3#4{\begin{figure}[tp]
	\begin{center}
		\includegraphics*[width=#1]{#2}
	\caption{#3}\label{#4} 
	\end{center}
	\end{figure}\par}

\begin{document}

    \input{sections}

    \bibliographystyle{ieee}
    \bibliography{literatu}

\end{document}

%% file: sections.tex
\title{\Large \bf Quadratic Control Lyapunov Functions for Bilinear Systems}

\author{Bernd Tibken, Frank Lehn and Eberhard P. Hofer}
\date{}
\maketitle
\thispagestyle{empty}

\vspace*{-8mm}
\begin{center}
Department of Measurement, Control and Microtechnology\\
University of Ulm, D-89069 Ulm, Germany\\
Phone: ++49 731 50-26300, Fax: ++49 731 50-26301
\end{center}

\vspace*{8mm}

\begin{abstract}
In this paper the existence of a quadratic control Lyapunov function
for bilinear systems is considered. The existence
of a control Lyapunov function ensures the existence of a control law which
ensures the global asymptotic stability of the closed loop control system. In this 
paper we will derive conditions for the existence of a control Lyapunov function
for bilinear systems. These conditions will be derived for the whole class of two dimensional bilinear systems with one control input. This will
lead to a simple flow diagram representation for the controller design for
bilinear systems using a quadratic control Lyapunov function. The controller design itself is carried out using Sontag's universal control law \cite{sontag:a} to obtain an asymptotically stable closed loop system. The Gutman control law \cite{gutman:a} which, however, only ensures practical stability of the closed loop system but which is considerably simpler than Sontag's control law will also be considered. An example will conclude the paper.
\end{abstract}

\section{Introduction}

In this paper bilinear systems (BLS) given by the state space representation
\begin{equation}\label{state}
\dot{x}(t)=Ax(t)+\sum_{i=1}^m(N_ix(t)+b_i)u_i(t)\,,\quad x(0)=x^0
\end{equation}
are in the focus of interest. In (\ref{state}) the symbol $x(t)\in R^n$ represents the state
vector and $u(t)\in R^m$ represents the vector of control inputs. 
The matrices $A$ and $N_i, i=1,\ldots,m$ are real
$n\times n$ matrices while the vectors $b_i, i=1,\ldots,m$ are real $n$-dimensional vectors. 
Applying the results from \cite{sontag:a, artstein:a, isidori:a} to the class of BLS
it turns out that $V(x)=x^TPx$ with a symmetric and positive definite $n \times n$ matrix $P$ 
is a control Lyapunov function (CLF) for (\ref{state}) if and only if 
\begin{equation}\label{condition}
Y(x):=x^T(A^TP+PA)x<0\quad \forall\,x\in M:=\left\{x\in R^n\left|\,\,x\neq 0\,\,,\,\,\sum_{i=1}^m((N_ix+b_i)^TPx)^2=0\right.\right\}
\end{equation} 
is valid. 
If $V(x)$ is a CLF the control law derived by Artstein and Sontag \cite{sontag:a, artstein:a} ensures the global asymptotic stability of the closed loop system. It should be mentioned that  
(\ref{condition}) is identical with the design conditions 
derived by Gutman \cite{gutman:a} and Mohler \cite{mohler:a} for the controller design for BLS.
However, in general the condition (\ref{condition}) is extremely hard to check and to fulfill.\\

In a numerical study, where a large set of bilinear systems has been considered we were in all cases unable to compute a matrix $P$ such that (\ref{condition}) is fulfilled. This numerical problem was the main 
motivation to analyse the design conditions for a special class of 
bilinear system, namely, two dimensional single input bilinear systems 
given by the state space representation
\begin{equation}\label{state-2D}
\dot{x}(t)=A\,x(t)+(Nx(t)+b)u(t),\,\, x(0)=x^0
\end{equation}
where the two dimensional state vector is represented by 
$x(t)=(x_1(t) , x_2(t))^T$ and the input $u(t)$ is a scalar function of time. 
The matrices $A$ and $N$ are real $2\times 2$ matrices and $b$ is a 
two dimensional vector. For this class of bilinear systems the computation of a positive definite 
and symmetric $2\times 2$ matrix $P=P^T>0$ for which
\begin{equation}\label{entwurfsbedingung}
Y(x)<0\quad\forall x\in M
\end{equation}
is satisfied will be carried out. In this special case the above definitions for $Y(x)$ and $M$ are simplified to
\begin{eqnarray}
Y(x)&=&x^T(A^TP+P\,A)x=x^TA_px\,,\label{yvonx}\\
M&=&\{x\mid x\neq 0,\quad (N\,x+b)^TP\,x=0\}=
\{x\mid x\neq 0,\quad x^TN_px+2x^TP\,b=0\}\label{mengem}
\end{eqnarray}
where the abbreviations
\begin{eqnarray}
A_p&:=&A^TP+P\,A\,,\label{apdef}\\
N_p&:=&N^TP+P\,N\label{npdef}
\end{eqnarray}
for the corresponding symmetric and real $2\times 2$ matrices 
have been introduced. The set $M$ defines a conic section in the 
$x$-plane whereby the shape of this conic section depends on the eigenvalues of 
$N_p$. For the analysis of the controller design a case study with respect 
to $N_p$ will be done. Especially positive definite, negative definite, semi definite and indefinite
matrices $N_p$ and the case $N_p=0$ have to be considered. 
In the following section a necessary condition for $P$ will be derived from (\ref{entwurfsbedingung}). 
The third section contains a detailed discussion of the positive definite
case $N_p>0$. The extension to the remaining cases will
be outlined and the results for all cases are presented as a flow diagram. The paper ends with an example and the conclusions. It should be mentioned that most of the computations were performed using the computer algebra system {\sl MAPLE} \cite{redfern:a}.

\section{Necessary Condition}\label{neccond}

The main goal in this section is to derive a necessary condition for $P$ 
from the design condition (\ref{entwurfsbedingung}). Because of $Y(0)=0$ the function $Y(x)$ has to 
have a local maximum at $x=0$ with respect to the constraint $x\in M$ because otherwise the 
condition (\ref{entwurfsbedingung}) is violated in a vincinity of $x=0$.
First we analyse if $Y(x)$ has a local extremum if the constraint 
$x^TN_px+2x^TP\,b=0$ which defines $M$ is active. To answer this, we define the Lagrange function 
$L(x,\lambda)=x^TA_px+\lambda(x^TN_px+2x^TPb)$ 
using the Lagrange multiplicator $\lambda$.
The necessary conditions for a local extremum using the 
partial derivatives are calculated as
\begin{equation}\label{lagequ}
\frac{\partial L}{\partial x}=2(A_p+\lambda N_p)x+2\lambda Pb\stackrel{!}{=}0,\quad
\frac{\partial L}{\partial \lambda}=x^TN_px+2x^TPb\stackrel{!}{=}0.
\end{equation}
It is easily verified that $\lambda=0$, $x=0$ is a solution of (\ref{lagequ}) and consequently $x=0$ is a local extremal point of $Y(x)$ under the constraint $x^TN_px+2x^TP\,b=0$. As discussed above 
$x=0$ is required to be a local maximum. A necessary and sufficient condition \cite{luenberger:a} for this is, that
\begin{equation}\label{hp1bed}
v^T(\frac{\partial^2}{\partial x^2}\{x^TA_px\})v<0
\end{equation}
is satisfied for all vectors $v$ that satisfy
\begin{equation}\label{hp2bed}
v^T(\frac{\partial}{\partial x}\{x^TN_px+2x^TP\,b\})\biggl|_{x=0}=0.
\end{equation}
The needed partial derivatives are $\frac{\partial^2}{\partial x^2}\{x^TA_px\}=2A_p$
and $\frac{\partial}{\partial x}\{x^TN_px+2x^TP\,b\}\bigl|_{x=0}=2P\,b$. The 
desired vectors $v$ that are othogonal to $P\,b$ and therefore satisfy (\ref{hp2bed}) 
are given by $v=\kappa J\,P\,b$ where $\kappa$ is a non zero scalar. In 
the two dimensional case the matrix $J$ is given by
\begin{equation}
J=\left(\begin{array}{cc}0&-1\\1&0  \end{array}\right).
\end{equation}
Now we are able to derive a condition that is equivalent to (\ref{hp1bed}) and 
(\ref{hp2bed}). It is given by
$v^TA_p v=\kappa^2b^TPJ^TA_pJPb\,<0\,\,\forall\kappa\in R\backslash\{0\}$ which can be simplified to
\begin{equation}\label{hproh}
b^TPJ^TA_pJPb<0
\end{equation}
since $\kappa^2$ is positive. This condition will be simplified further. In the following we assume the pair $\{A,b\}$ to be completely controllable and choose $\{A,b\}$ in controller normal form such that
\begin{equation}\label{abnormal}
A=\left(\begin{array}{cc}0&1\\-a_0&-a_1\end{array}\right),\,\,\,\,
b=\left(\begin{array}{c}0\\1\end{array}\right)
\end{equation}
where $a_0$ and $a_1$ are the coefficients of the characteristic polynomial 
of $A$. The positive definite and symmetric matrix $P$ is parameterized as
\begin{equation}\label{pparam}
P=\left(\begin{array}{cc}1&p_1\\p_1&p_2\end{array}\right)
\end{equation}
where $p_1$ and $p_2$ are real variables and the fact that $P$ is determined 
up to a positive factor has been used to choose the $(1,1)$ element to 1. To 
ensure $P>0$ the variables $p_1$ and $p_2$ have 
to fulfill
\begin{equation}\label{pposdef}
p_2-p_1^2>0.
\end{equation}
If we use (\ref{abnormal}) and (\ref{pparam}) in (\ref{hproh}) the condition 
(\ref{hproh}) is equivalent to $p_1(p_1^2-p_2)<0$. Using (\ref{pposdef}) this 
inequality reduces to
\begin{equation}\label{p1gl}
p_1>0
\end{equation}
which is a necessary condition on $P$ and will be central for the following computations.

\section{Sufficient Conditions}\label{scond}

In this section we discuss the controller design conditions for the case of
a positive definite matrix $N_p$. The main idea to analyse this situation is 
to carry out a linear state variable transformation to transform the set $M$ 
into a circle. After this transformation we pa\-ra\-me\-te\-ri\-ze $M$ and evaluate 
$Y(x)\,\,\forall x\in M$. In the following the transformed set $M$ and the 
transformed function $Y(x)$ will be called $M_y$ and $Y_y(y)$, respectively. 
It should be mentioned that temporarily $x=0$ will be allowed in $M$ for an
easier calculation of the necessary transformation. After having carried 
out the transformation it will be indicated to which point $x=0$ has been transformed and the 
related transformed point will be omitted. For the positive definite 
matrix $N_p$ we carry out a cholesky decomposition \cite{lancaster:a} of the 
form $N_p=L^TL$ with a nonsingular matrix $L$ given by
\begin{equation}
L=\left(\begin{array}{cc}l_1&l_2\\0&l_3\end{array}\right)
\end{equation}
where $l_1$ and $l_3$ are positive. With this matrix $L$ the above mentioned 
transformation is performed according to
$x=L^{-1}(y+y_0)\Leftrightarrow y=Lx-y_0$
where $y=(y_1,y_2)^T$ represents the new state vector and
\begin{equation}
y_0=-(L^{-1})^TPb=\left(\begin{array}{c} 
-\frac{p_1}{l_1}\\ \frac{p_1l_2-p_2l_1}{l_1l_3}    
\end{array}\right)
=\left(\begin{array}{c}
y_{0,1}\\y_{0,2}
\end{array}\right)
\end{equation}
is an appropriate choice. Substituting $x=L^{-1}(y+y_0)$ into the equations 
which describe $M$ leads to the transformed set
\begin{equation}
M_y=\{y\in R^2\mid y_1^2+y_2^2=a\}
\quad\mbox{with}\quad
a=\frac{(p_1l_2-p_2l_1)^2+(p_1l_3)^2}{l_1^2l_3^2}>0
\end{equation}
which describes a circle with radius $\sqrt{a}$ in the $y$-plane. This circle 
is now parameterized using rational functions, thus, the transformed set $M_y$ 
is represented by
\begin{equation}
M_y=\{y\mid y=\frac{\sqrt{a}}{1+t^2}\left(\begin{array}{c}1-t^2\\\ 2t\end{array}\right), t\in R\}.
\end{equation}
After having computed this representation of $M_y$ we are now able to 
evaluate $Y(x)\,\,\forall x\in M$, which has the same range of values as 
the transformed
$Y_y(y)\,\,\forall y\in M_y$. We compute $Y_y(y)=Y(L^{-1}y+x_0)$
and in order to evaluate $Y_y(y)\,\,\forall y\in M_y$ we further substitute
\begin{equation}
y=\frac{\sqrt{a}}{1+t^2}\left(\begin{array}{c}1-t^2\\\ 2t\end{array}\right)
\end{equation}
in $Y_y(y)$ and compute $Y_y(y)\,\forall y\in M_y$ as a rational 
function $\tilde{Y}_y(t)$ of the real variable $t$. The function 
$\tilde{Y}_y(t)$ is required to be negative for all real values of the 
variable $t$ to ensure that the design condition (\ref{entwurfsbedingung}) holds. $\tilde{Y}_y(t)$ is a very complex rational function of the form 
$\tilde{Y}_y(t)=\frac{Z(t)}{N(t)}$, where numerator 
$Z(t)$ and denominator $N(t)$ are both polynomials of degree four in $t$.
The denominator polynomial is computed as $N(t)=l_1^4l_3^4(1+t^2)^2$ which is 
positive for all real values of $t$, thus, it can be ignored for the  
analysis of the sign of $\tilde{Y}_y(t)$. The numerator polynomial $Z(t)$ is 
factored as the product of two polynomials of degree two given by
\begin{equation}\label{t0eq}
Z(t)=(t-t_0)^2\cdot P(t)
\quad\mbox{with}\quad
t_0=-\frac{-y_{0,2}}{\sqrt{a}+(-y_{0,1})}.
\end{equation}
The special value $t_0$ corresponds to $x=0$ and because $x=0$ has to be excluded from $M$ we require $Z(t)<0\,\,
\forall\,t\in R\backslash\{t_0\}$. For all values of $t\in R\backslash
\{t_0\}$ the polynomial $(t-t_0)^2$ is positive so it does not influence the 
sign of $\tilde{Y}_y(t)$. Thus, $P(t)$ is required to be negative for all 
real values of $t$ in order to satisfy the design condition (\ref{entwurfsbedingung}). As mentioned 
above $P(t)$ is a polynomial of degree two of the form $P(t)=k_2t^2+k_1t+k_0$, 
where $k_2$, $k_1$ and $k_0$ are very complex expressions containing $l_1$, 
$l_2$, $l_3$, $p_1$, $p_2$, $a_0$ and $a_1$. If we require $P(t)<0$ for all 
$t\in R$ the following two cases have to be considered:
\begin{enumerate}
\item 	all coefficients $k_2$, $k_1$ and $k_0$ are non zero, thus,
	\begin{equation}\label{normpoly}
		P(t^*)<0\,\,\mbox{for some}\,\,t^*\,\,\mbox{and}\,\,
		k_1^2-4k_0k_2<0
	\end{equation}
	have to be fulfilled, in order to guarantee  $P(t)<0$ 
	for all real $t$;
\item 	the coefficients $k_2$ and $k_1$ are zero and therefore $P(t)=k_0$ is 
	constant and $k_0<0$ has to be fulfilled.
\end{enumerate}
These two cases will be discussed in the following. First case one is 
considered. If $t^*=t_0$ given by (\ref{t0eq}) is chosen, the conditions (\ref{normpoly}) are computed as
\begin{eqnarray}
P(t_0)=8l_1^2l_3^2p_1(p_1^2-p_2)&<&0,\label{pvt0}\\
k_1^2-4k_0k_2=16\beta l_1^2l_3^2(\sqrt{\beta}+p_1l_3)^2\cdot&&\nonumber\\
\cdot(4p_1^2a_0+p_1^2a_1^2-2p_1p_2a_0a_1-
2p_1a_1-2p_2a_0+p_2^2a_0^2+1)&<&0\label{diskrim}
\end{eqnarray}
with $\beta=(p_1l_2-p_2l_1)^2+(p_1l_3)^2$. Considering the necessary conditions (\ref{pposdef}) 
and (\ref{p1gl}) one can easily conclude that condition (\ref{pvt0}) is satisfied. 
It is also obvious, that the factor $16\beta l_1^2l_3^2(\sqrt{\beta}+p_1l_3)^2$ 
in (\ref{diskrim}) is positive and consequently, 
\begin{equation}\label{gkr_roh}
(4p_1^2a_0+p_1^2a_1^2-2p_1p_2a_0a_1-
2p_1a_1-2p_2a_0+p_2^2a_0^2+1)<0
\end{equation}
is the only condition that is left. Condition (\ref{gkr_roh}) is now rewritten 
in two different forms, which are given by
\begin{eqnarray}
\left[(a_1p_1-a_0p_2)-1\right]^2+4a_0(p_1^2-p_2)&<&0\label{sw1}\\
\mbox{and}&&\nonumber\\
\left[(a_0p_2-a_1p_1)-1\right]^2+4p_1(a_0p_1-a_1)&<&0\label{sw2}.
\end{eqnarray}
From (\ref{sw1}) it is concluded that $a_0>0$ has to hold if we 
consider $p_2-p_1^2>0$ from (\ref{pposdef}). If we now use $a_0>0$ in 
(\ref{sw2}) it is obvious that $a_1>0$ has to hold if we require 
$p_1>0$ from (\ref{p1gl}). The following theorem has therefore been proved.\\

{\bf Theorem} {\it The condition
$(4p_1^2a_0+p_1^2a_1^2-2p_1p_2a_0a_1-
2p_1a_1-2p_2a_0+p_2^2a_0^2+1)<0$
with respect to the conditions (\ref{pposdef}) and (\ref{p1gl}) can be fulfilled if and 
only if $a_0>0$ and $a_1>0$ hold. Thus, the bilinear system with $u\equiv0$ 
has to be asymptotically stable in order to allow the controller design in this case.}\\

In the following we consider case two ($k_2=k_1=0$ and consequently $k_0<0$).
Since $k_2$ and $k_1$ are both linear functions of $a_0$ and $a_1$ we use 
$k_2=0$ and $k_1=0$ as a linear set of equations to calculate $a_0$ and $a_1$. 
In other words: For which bilinear system (defined by $a_0$ and $a_1$) is it possible, that this case
occurs. The solution with respect to $a_0$ and $a_1$ is given by
\begin{eqnarray}
a_0&=&\frac{1}{\sqrt{\beta}N_{k_0}}\left(p_2-p_1^2\right)
\left(l_2p_1-p_2l_1\right)^2l_1^2,\label{a0gl}\\
a_1&=&\frac{1}{\beta N_{k_0}} g(l_1,l_2,l_3,p_1,p_2)\label{a1gl}
\end{eqnarray}
where $g(l_1,l_2,l_3,p_1,p_2)$ and $N_{k_0}$ are very complicated expressions depending on $l_1$, $l_2$, $l_3$, $p_1$ and $p_2$.
These values for $a_0$ and $a_1$ are substituted into the expression for $k_0$ 
and we compute
\begin{equation}
k_0=\,-\,\frac{8\sqrt{\beta}p_1l_1^2l_3^2}{N_{k_0}}\left(\sqrt{\beta}+
p_1l_3\right)^2\left(p_2-p_1^2\right).
\end{equation}
The only interesting case appears if $N_{k_0}>0$ holds since for this case 
$k_0$ is negative because the necessary conditions (\ref{pposdef}) and (\ref{p1gl}) have to be 
fulfilled. Considering this case ($N_{k_0}>0$), one can easily see from (\ref{a0gl}) that 
$a_0>0$ is valid if $l_2p_1-p_2l_1\neq 0$ holds. 
In the following the sign of $a_1$ from (\ref{a1gl}) is analyzed for the case $l_2p_1-p_2l_1\neq 0$, which reduces to the analysis of $g(l_1,l_2,l_3,p_1,p_2)$ since $\beta>0$ holds and $N_{k_0}>0$ has turned out to be the only interesting case. The case $l_2p_1-p_2l_1=0$ will be discussed later. 
Introducing $\mu_1=p_1l_2$, $\mu_2=p_2l_1>0$ 
and $\mu_3=p_1l_3>0$ and using polar coordinates $\mu_3=r\cos(\varphi)$ and 
$\mu_1=\mu_2+r\sin(\varphi)$ with $r>0$ and $-\frac{\pi}{2}<\varphi<\frac{\pi}{2}$ 
we compute
\begin{equation}
\beta =(p_1l_2-p_2l_1)^2+(p_1l_3)^2=(\mu_1-\mu_2)^2+\mu_3^2=r^2 \quad
\Leftrightarrow\quad \sqrt{\beta}=r
\end{equation}
and $g(l_1,l_2,l_3,p_1,p_2)=\tilde{g}(r,\varphi,\mu_2)$ is given by
$\tilde{g}(r,\varphi,\mu_2)=4r^4\cos(\varphi)(1+\cos(\varphi))\,\,
\mbox{$q^T Q q$}$. Hereby, the two dimensional vector $\mbox{$q$}=\left(\begin{array}{cc}
r&\mu_2 \end{array}\right)^T$ and the $2\times 2$ matrix
\begin{equation}
\mbox{$Q$}=\left(\begin{array}{cc}
q_1	&q_2	\\
q_2 	&q_3	
\end{array}\right)
\end{equation}
with
\begin{equation}
q_1=p_2\cos^2\left(\frac{\varphi}{2} \right),\quad
q_2=(p_2-p_1^2)\cos\left(\frac{\varphi}{2} \right)\sin\left(\frac{\varphi}{2} \right),\quad
q_3=(p_2-p_1^2)\sin^2\left(\frac{\varphi}{2} \right)
\end{equation}
which is positive definite if $\varphi\neq 0$ have been introduced. Looking at this representation of $\tilde{g}(r,\varphi,\mu_2)$ as a 
quadratic form with positiv definite matrix $\mbox{$Q$}$ it is 
obvious that $g(l_1,l_2,l_3,p_1,p_2)>0$ and consequently $a_1>0$ from 
(\ref{a1gl}) holds for $\varphi\neq 0 \Leftrightarrow l_2p_1-p_2l_1\neq0$. 
The case $l_2p_1-l_1p_2=0$ which belongs to $\varphi=0$ will be discussed in 
the following. 
For $\varphi=0$ the matrix $\mbox{$Q$}$
is only positive semidefinite and the previous conclusions are no longer valid. 
Using $l_2p_1-l_1p_2=0$ and substituting $l_2=\frac{p_2}{p_1}l_1$ into (\ref{a0gl}) and (\ref{a1gl}) leads to
\begin{equation}
a_0=0,\quad
a_1=\frac{1}{p_1}
\end{equation}
which corresponds to a matrix $A$ that is not asymptotically stable. 
Resuming the case of a positive definite matrix $N_p$, the controller design is only possible if $a_0>0$ and
$a_1>0$ hold if $l_2p_1-p_2l_1\neq 0$ is valid.
For the case $p_1l_2-p_2l_1 = 0$ a system with 
$a_0=0$ and $a_1=\frac{1}{p_1}>0$ can be allowed. In this case we arrived 
at an interesting point since we are able to design a controller for a 
bilinear system which is not asymptotically stable for $u\equiv 0$. The remaining cases for $N_p$
are treated in the same style, namely, after the transformation of $M$ into a 
representative set as circle, parabola or hyperbola the set is parameterized
and $Y$ is investigated. In the next section a flow diagram for controller
design is presented containig all possible cases for $N_p$.

\section{Flow Diagram}\label{flowdiagram_sect}
In this section we give a flow diagram
for the controller design for a two dimensional bilinear system
with one input $u$. For the design we assume a system  that is given in 
controller normal form such that the matrices $A$ and $b$ have the structure 
defined in (\ref{abnormal}). Each system where the pair $\{A,b\}$ is completely controllable
can easily be transformed into this normal form using a linear state space 
transformation. The results of the preceeding section and of the other cases to be treated
are summarized in the flow diagram shown in Figure \ref{fdiag}.
Using this flow diagram it is possible to calculate the required matrix $P$ and carry out
the controller design according to, e.g. Sontag \cite{sontag:a} or Gutman \cite{gutman:a}.

\fpicps{15cm}{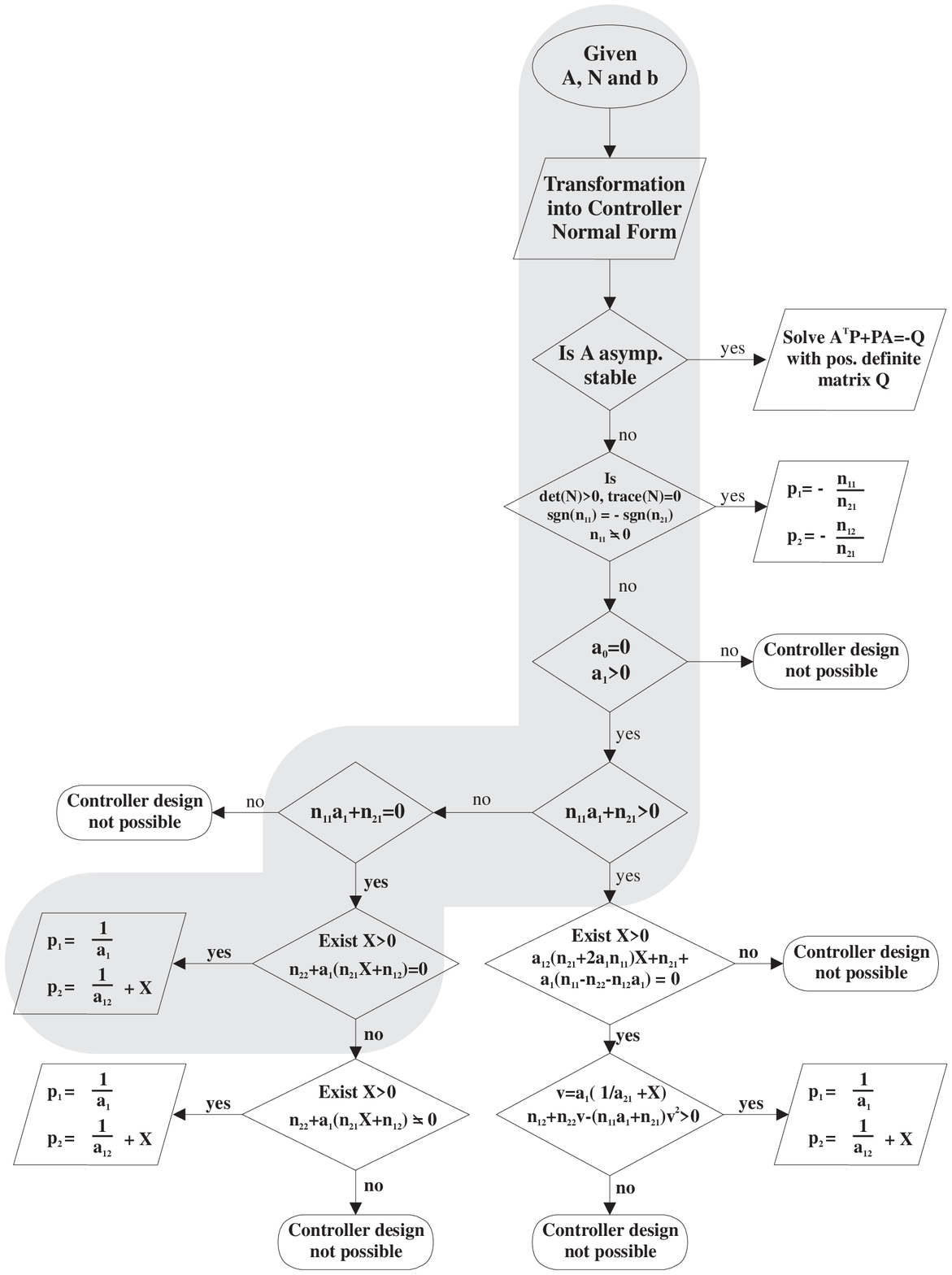}{Flow Diagram for the Controller Design.}{fdiag}

\section{Example}\label{example}
We consider the 
two dimensional bilinear system in controller normal form with the system matrices
\begin{equation}\label{testsystem}
A=\left(\begin{array}{cc}
0&1\\0&-1
\end{array}\right),\,\,\,
N=\left(\begin{array}{cc}
1&1\\-1&1
\end{array}\right)\,\,\mbox{and}\,\,
b=\left(\begin{array}{c}
0\\1
\end{array}\right).
\end{equation}
Now we step through the flow diagram and have to answer the following 
questions which correspond to the shaded path in Figure \ref{fdiag}.
\begin{enumerate}
\item Is the system asymptotically stable?  We compute the characteristic polynomial 
$p(s)=s(s+1)$ and the eigenvalues of $A$ as $s_1=0$ and $s_2=-1$; the system is not asymtotically stable.
\item $\mbox{det}(N)>0$, $\mbox{trace}(N)=0$, $n_{11}\neq 0$ and 
$\mbox{sgn}(n_{11})=-\mbox{sgn}(n_{21})$?  We compute $\mbox{det}(N)=2>0$,
$\mbox{trace}(N)=2$, $n_{11}=1\neq 0$, $\mbox{sgn}(n_{11})=1$, and
$\mbox{sgn}(n_{21})=-1$; these conditions are not fulfilled, since $\mbox{trace}(N)\neq 0$.
\item Is $a_0=0$ and $a_1>0$? We consider the matrix $A$ and read $a_0=0$
and $a_1=1>0$; these conditions hold.
\item Is $n_{11}a_1+n_{21}>0$? We compute $n_{11}a_1+n_{21}=0$; we continue with the "no" branch.
\item Is $n_{11}a_1+n_{21}=0$? We compute $n_{11}a_1+n_{21}=0$; this condition holds.
\item Exist $X>0$ such that $n_{22}+a_1(n_{21}X+n_{12})=0$? We evaluate
this equation and compute $X=2>0$; we are able to compute the desired
positive value of $X$.
\item Compute $p_1$ and $p_2$! We compute
$p_1=\frac{1}{a_1}=\frac{1}{1}=1$ and 
$p_2=\frac{1}{a_1^2}+X=\frac{1}{1}+2=3$. 
\item Set up the matrix $P$! We compute
\[
P=\left(\begin{array}{cc}
1&p_1\\p_1&p_2
\end{array}\right)\,\,=\,\,
\left(\begin{array}{cc}
1&1\\1&3
\end{array}\right).
\]
\item Compute the control law. In this example the control law given by Gutman \cite{gutman:a} 
\begin{equation}\label{gutmancontrollaw}
u=-\alpha(Nx+b)^TPx=-\alpha(4x_2^2+x_1+3x_2)
\end{equation}
is used.
\end{enumerate}
Now we will check if the condition (\ref{entwurfsbedingung}) holds. Therefore we have to evaluate
$Y=x^T(A^TP+PA)x$ for all vectors $x$ that belong to the set
$M=\{x\in R^2\,\mid\, x\neq 0,\,x^T(N^TP+PN)x+2x^TPb=0\}$. 
We first compute the set $M$ using $x=(x_1\,,\,x_2)^T$ and we have
$M=\{x\in R^2\,\mid\, x\neq 0,\, 4x_2^2+x_1+3x_2=0\}$.
Now we evaluate $Y$ for all $x\in M$ and compute $Y\mid_{x\in M}=-4x_2^2,\,\, x_2\neq 0$.
This function is negative for all real values of $x_2\neq 0$. Thus, using this matrix $P$ the design condition (\ref{entwurfsbedingung}) holds and $V(x)=x^TPx$ is a CLF for the BLS (\ref{testsystem}). Figure \ref{sim} 
shows trajectories of the closed loop system for the choice $\alpha=\frac{1}{10}$ in (\ref{gutmancontrollaw}).
These trajectories clearly indicate the global asymptotic stability in this case.
\fpicps{11.5cm}{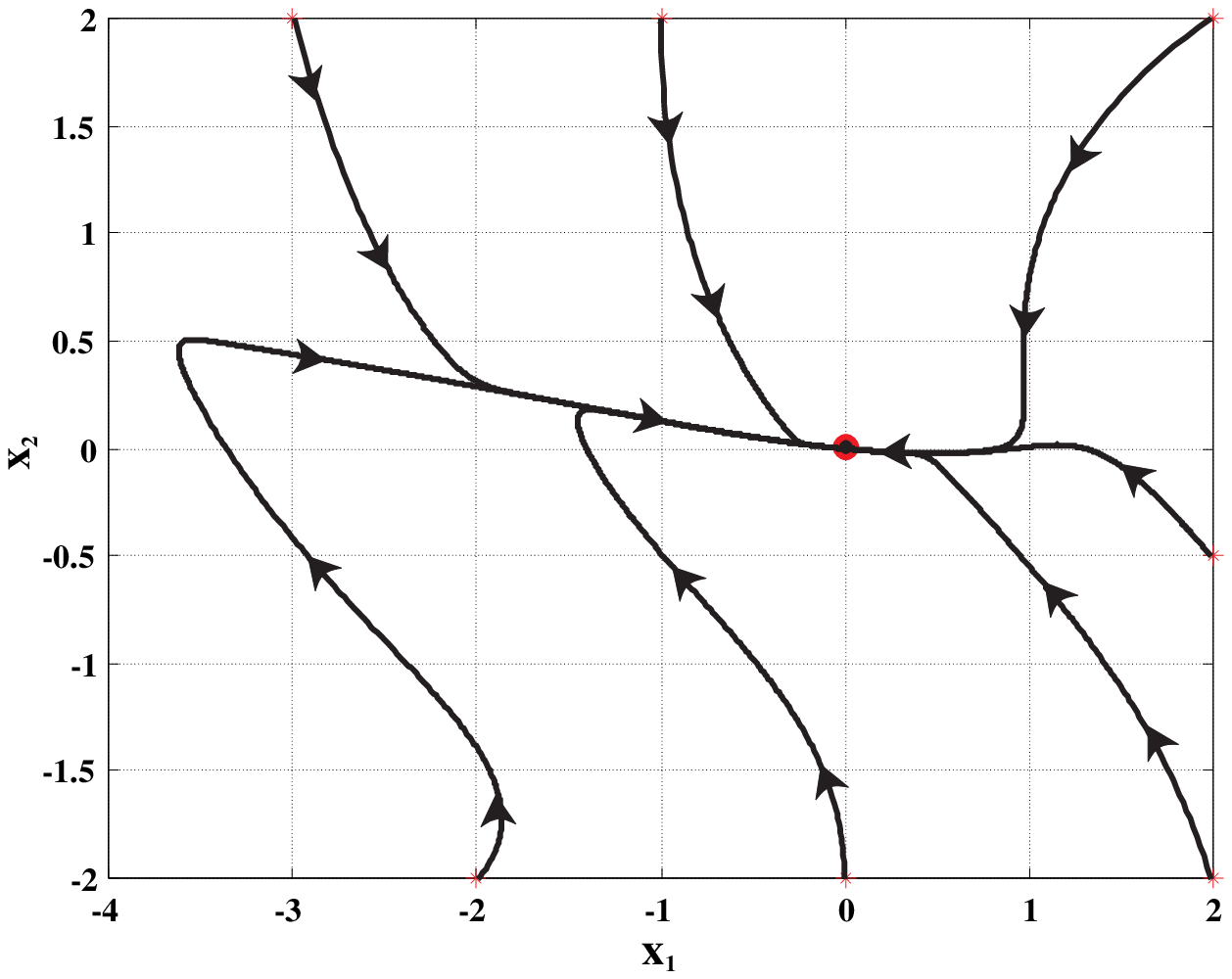}{Simulation result for the closed loop system.}{sim}

\section{Conclusions}\label{conc}
The controller design for bilinear systems is a very complicated task 
and only a few results for higher dimensional systems are available.
Most of the recent publications treat systems with asymptotically 
stable system matrix in the linear part.
In this paper the existence of a quadratic CLF for bilinear systems has been
analysed systematically for the case of two dimensional systems
with one input. This case is important because no systematic
study for two dimensional systems has been carried out so far
and the conditions (\ref{condition}) simplify with increasing
number of inputs. Therefore, two dimensional systems with one
input represent the first interesting case which has not been
treated systematically up to now. We give
necessary and sufficient conditions 
which ensure that the design condition (\ref{entwurfsbedingung}) can be
satisfied. The result of the
complete analysis is presented as a flow diagram which guides
the user through the controller design. This is the first
complete systematic study of the conditions (\ref{condition}) for a
whole class of bilinear systems. 
In previous research only
special examples have been treated or the case of bilinear
systems with asymptotically stable system matrix of the linear
part has been investigated. 
The important case with an unstable
system matrix of the linear part has not been treated in the
literature so far. This gap is closed by the results of this
paper for the case of two dimensional bilinear systems with one
input. An example illustrating the design is given at the end of
the paper.